\documentclass[11pt]{amsart}

\usepackage{amsmath}
\usepackage{amssymb}
\usepackage{mathrsfs}
\usepackage[all]{xy}
\usepackage[utf8]{inputenc}

\newtheorem{theorem}{Theorem}[section]
\newtheorem{claim}[theorem]{Claim}

\newtheorem{lemma}[theorem]{Lemma}

\newtheorem{proposition}[theorem]{Proposition}
\newtheorem{corollary}[theorem]{Corollary}

\theoremstyle{definition}
\newtheorem{definition}[theorem]{Definition}

\newtheorem{question}[theorem]{Question}

\theoremstyle{remark}
\newtheorem{remark}[theorem]{Remark}

\newcount\skewfactor
\def\mathunderaccent#1#2 {\let\theaccent#1\skewfactor#2
\mathpalette\putaccentunder}
\def\putaccentunder#1#2{\oalign{$#1#2$\crcr\hidewidth
\vbox to.2ex{\hbox{$#1\skew\skewfactor\theaccent{}$}\vss}\hidewidth}}
\def\name{\mathunderaccent\tilde-3 }

\def\smallbox#1{\leavevmode\thinspace\hbox{\vrule\vtop{\vbox
   {\hrule\kern1pt\hbox{\vphantom{\tt/}\thinspace{\tt#1}\thinspace}}
   \kern1pt\hrule}\vrule}\thinspace}

\newcommand{\cf}{{\rm cf}}

\newcommand{\stick}{\ensuremath\mspace{2mu}\mid \mspace{-12mu}{\raise 0.4em \hbox{$\bullet$}}}

\def\qedref#1{$\qed_{\reforiginal{#1}}$}

\title{Stationary and Closed Rainbow subsets}
\author{Shimon Garti}
\address{Einstein Institute of Mathematics,
 The Hebrew University of Jerusalem,
 Jerusalem 91904, Israel}
\email{shimon.garty@mail.huji.ac.il}

\author{Jing Zhang}
\address{Department of Mathematics, Bar-Ilan University, 5290002, Israel}
\email{jingzhangjz13@gmail.com}

\thanks{The second author is supported by the Foreign Postdoctoral Fellowship Program of the Israel Academy of Sciences and Humanities and by the Israel Science Foundation (grant agreement 2066/18).}


\begin{document}
\let\labeloriginal\label
\let\reforiginal\ref

\begin{abstract}
We study the structured rainbow Ramsey theory at uncountable cardinals. When compared to the usual rainbow Ramsey theory, the variation focuses on finding a rainbow subset that not only is of a certain cardinality but also satisfies certain structural constraints, such as being stationary or closed in its supremum. In the process of dealing with cardinals greater than $\omega_1$, we uncover some connections between versions of Chang's Conjectures and instances of rainbow Ramsey partition relations, addressing a question raised in \cite{zhang}.
\end{abstract}

\maketitle

\section{Introduction}

The study of rainbow Ramsey theory is dual to that of the usual Ramsey theory. A typical problem is that: given $f:[\kappa]^n\rightarrow\theta$ satisfying certain constraints, we are asked to find $y\subseteq\kappa$ satisfying certain constraints, such that $y$ is $f$-\emph{rainbow}, namely, $f\upharpoonright[y]^n$ is one-to-one.

\begin{definition}\label{rainbowinitialdef}
We use $\lambda\to^{\rm{poly}} (\kappa)^n_{<i-bdd}$ to abbreviate: for any $f: [\lambda]^n \to \lambda$ that is $<i$\emph{-bounded}, namely for any $\alpha\in \lambda$, $|f^{-1}\{\alpha\}| < i$, there exists an $f$-rainbow $A\subset \lambda$ of order type $\kappa$.
\end{definition}

We will let $\lambda\to^{\rm{poly}} (\kappa)^n_{i-bdd}$ abbreviate $\lambda\to^{\rm{poly}} (\kappa)^n_{<i^+-bdd}$.

Initial results in the area of infinite rainbow Ramsey theory were obtained by Galvin, followed by a series of results that appeared in \cite{MR716846}, \cite{MR2354904}, \cite{MR2902230}, \cite{zhang}. For a short history and an account of basic facts, we direct the reader to  \cite{MR2354904}.

In this paper, we are concerned with structural strengthenings of the usual rainbow Ramsey partition relations, loosely motivated by the study of \emph{topological partition relations}, see \cite{MR2912439} for example.

Boundedness conditions based on the cardinality as in Definition \ref{rainbowinitialdef} are the most studied constraints for the given coloring. We will encounter other types of boundedness conditions, as an immediate step to solve problems involving the cardinality-based boundedness conditions. The following general definition captures the main problems we would like to work on.
 Recall that a set $A\subset \rm{Ord}$ is \emph{closed}, if whenever $\sup A\cap \alpha = \alpha<\sup A$, then $\alpha\in A$.

\begin{definition}\label{generalproperty}
Let $\mathcal{P}$ be a property for colorings on $[\lambda]^n$. 
\begin{enumerate}
\item We use $\lambda\to^{\rm{poly}} (\alpha-cl)^n_\mathcal{P}$ to abbreviate: for any $f: [\lambda]^n\to \lambda$ satisfying $\mathcal{P}$, there exists a \emph{closed} $f$-rainbow subset of order type $\alpha$. 
\item Similarly, $\lambda\to^{\rm{poly}} (\alpha-st)^n_\mathcal{P}$ abbreviates: for any $f: [\lambda]^n\to \lambda$ satisfying $\mathcal{P}$, there exists a $f$-rainbow subset $A$ of order type $\alpha$ that is \emph{stationary in $\sup A$}. 
\end{enumerate}
\end{definition}

Clearly, $\lambda\to^{\rm{poly}} (\alpha-st)^n_\mathcal{P}$ only makes sense if $\cf(\alpha)>\omega$. Todor\v{c}evi\'{c} in  \cite{MR716846} showed that it is consistent relative to the consistency of \textsf{ZFC} that $\omega_1 \to^{\rm{poly}} (\omega_1-st)^2_{<\omega-bdd}$ and Abraham, Cummings and Smyth in \cite{MR2354904} showed that under strong forcing axioms, namely Martin's Maximum, $\omega_2 \to^{\rm{poly}} (\omega_1-cl)^2_{<\omega-bdd}$. We will consider the improvements where the boundedness conditions are relaxed. Clearly, $\omega_1\not\to^{\rm{poly}} (\omega_1-st)^2_{\omega-bdd}$ outright in \textsf{ZFC}. We show the next natural instance, after enlarging the source cardinal, is consistent.

Our main result reads: 
\begin{theorem}\label{mainstat}
Relative to the existence of a huge cardinal, it is consistent that $\omega_2\to^{\rm{poly}} (\omega_1-st)^2_{\omega-bdd}$.
\end{theorem}

The consistency of $\omega_2\to^{\rm{poly}} (\omega_1-cl)^2_{\omega-bdd}$ relative to any reasonable assumption is still open. However, we are able to demonstrate an improvement to the Abraham-Cummings-Smyth result by relaxing the boundedness condition.

An interesting role that versions of Chang's Conjectures plays in certain ``stepping-up'' arguments has also been uncovered. Let us recall some definitions. Suppose that $\theta,\kappa,\lambda,\mu$ are infinite cardinals. \emph{Chang's Conjecture} $(\mu,\lambda)\twoheadrightarrow(\kappa,\theta)$ is the following assertion: for every countable first-order language $\mathcal{L}$ and for every $\mathcal{L}$-structure $A=(\mu, \lambda, .... )$ there exists an $X\prec A$ such that $|X|=\kappa$ and $|X\cap \lambda|=\theta$.
The special case $(\omega_2,\omega_1)\twoheadrightarrow(\omega_1,\omega)$ is called simply the \emph{Chang's Conjecture} (CC).

In \cite{MR2354904} and \cite{MR2902230}, it was shown that if $\kappa$ is a regular cardinal satisfying $\kappa^{<\kappa}=\kappa$, then $\kappa^+\to^{\rm{poly}} (\eta)^2_{<\kappa-bdd}$ for all $\eta<\kappa^+$. Furthermore, these positive relations persist in certain forcing extensions. A question regarding the consistency of the partition relations when $\kappa$ is singular was asked. In \cite{zhang}, the following were shown: for any singular $\kappa$ with $\lambda=\cf(\kappa)<\kappa$,
\begin{enumerate}
\item $\kappa^+\not\to^{\rm{poly}} (\eta)^2_{<\kappa-bdd}$, whenever $\eta>\lambda^+$, 
\item $\kappa^+\to^{\rm{poly}} (\eta)^2_{<\kappa-bdd}$ for all $\eta<\lambda^+$ if in addition $\kappa^{<\lambda}=\kappa$ and
\item $S^{\kappa^+}_{\lambda^+} \in  I[\kappa^+]$ implies $\kappa^+\not\to^{\rm{poly}} (\lambda^+)^2_{<\kappa-bdd}$.
\end{enumerate}

Here $I[\kappa^+]$ is Shelah's \emph{approachability ideal}. See \cite{MR2768694} for more information on these ideals. In particular, (3) above shows that $\kappa^+\to^{\rm{poly}} (\lambda^+)^2_{<\kappa-bdd}$ has high large cardinal strength, if consistent at all. A question regarding the consistency of $\kappa^+\to^{\rm{poly}} (\omega_1)^2_{<\kappa-bdd}$ for some singular $\kappa$ of countable cofinality was asked in \cite[Question 3.16]{zhang}. We will answer this question positively utilizing a version of the singular Chang's Conjectures whose consistency was first established in \cite{MR1045371}.

The structure of the paper is: 

\begin{enumerate}
\item Section \ref{warmup} contains some general discussion on obtaining closed rainbow subsets of various order types,
\item Section \ref{CC} contains the proof for Theorem \ref{mainstat},
\item Section \ref{pfanotimplies} demonstrates that PFA does not imply $\omega_2\to^{\rm{poly}} (\omega_1-st)^2_{\omega-bdd}$,
\item Section \ref{mmimplies} illustrates a strengthening of $\omega_2 \to^{\rm{poly}} (\omega_1-cl)^2_{<\omega-bdd}$ under \textsf{MM}, 
\item Section \ref{openquestion} concludes with some open questions.
\end{enumerate}

Unless otherwise stated, we will assume all the colorings are \emph{normal}, in the sense that a coloring $f$ whose domain is $[\kappa]^n$ is normal, if whenever $a, b\in [\kappa]^n$ satisfy that $f(a)=f(b)$, we have $\max a = \max b$.

\section{Closed rainbow sets}\label{warmup}

We begin with discussing some \textsf{ZFC} constraints. 
\begin{proposition}\label{constraint}
For any cardinal $\kappa$ with $\cf(\kappa)>\omega$, $\kappa\not\to^{\rm{poly}}(\cf(\kappa)-cl)^2_{2-bdd}$.
\end{proposition}

Recall that it is a theorem of Shelah \cite[Claim 2.3]{MR1318912} that if $\kappa=\cf(\kappa)\geq \omega_2$ then there is a club-guessing sequence supported on $S^\kappa_\omega$. The consequence that we will use is that there exists a sequence $\langle C_\alpha : \alpha\in S^\kappa_\omega\rangle$ where each $C_\alpha$ of order type $\omega$ is cofinal in $\alpha$, so that any club $D\subset \kappa$, there exists some $\alpha\in S^\kappa_\omega$ such that $C_\alpha\subset D$.

\begin{lemma}
\label{thmomega2closed}
Suppose that $\kappa=\cf(\kappa)\geq\omega_2$. Then $\kappa\nrightarrow^{\rm poly}(\kappa-cl)^2_{2-bdd}$.
\end{lemma}

\par\noindent\emph{Proof}. \newline 
Let $S=S^{\kappa}_\omega$ and let $\bar{A}=(A_\delta:\delta\in S)$ be a club guessing sequence.
We define a normal 2-bounded coloring $f$ by induction on $\delta\in\kappa$, making sure that for every $\delta\in S$, there are $\alpha_\delta \neq \beta_\delta\in A_\delta$ satisfying $f(\alpha_\delta, \delta)=f(\beta_\delta, \delta)$.

Let $E\subseteq\kappa$ be a club and let $S_E=\{\delta\in S:A_\delta\subseteq E\}$. Since $\bar{A}$ guesses clubs, $S_E$ is stationary.
Choose $\delta\in \lim (E\cap S_E)\cap (E\cap S_E)$.
By the construction we have $\alpha_\delta \neq\beta_\delta\in A_\delta\subseteq E$ so that $f(\alpha_\delta,\delta)=f(\beta_\delta,\delta)$.
Since $\alpha_\delta,\beta_\delta,\delta\in E$, we conclude that $E$ is not $f$-rainbow.

\hfill \qedref{thmomega2closed}

If there is a club guessing sequence $(A_\delta:\delta\in S^{\omega_1}_\omega)$ then $\omega_1\nrightarrow^{\rm poly}(\omega_1-cl)^2_{2-bdd}$ by the same reasoning. However, \textsf{ZFC} does not prove the existence of such club guessing sequences. Hence, we employ a different proof in the case of a successor cardinal.

%Now if $\kappa=\cf(\kappa)$ and $\lambda=\cf(\lambda)>\kappa^+$ then we have a club guessing in \textsf{ZFC}, but if $\lambda=\kappa^+$ then one can force that there is no club guessing for $\kappa$ and $\lambda$.
%In particular, under \textsf{MM} (or just \textsf{PFA}) there is no club guessing on $S^{\omega_1}_\omega$.
%Nevertheless, we can prove the following:

\begin{lemma}
\label{clmomega1} If $\kappa$ is an infinite cardinal then $\kappa^+\nrightarrow^{\rm poly}(\kappa^+-cl)^2_{2-bdd}$.
\end{lemma}

\par\noindent\emph{Proof}. \newline 
Suppose for the sake of contradiction that $\kappa^+\rightarrow^{\rm poly}(\kappa^+-cl)^2_{2-bdd}$.
For each ordinal $\beta\in\kappa^+ - \kappa$, fix a bijection $f_\beta:\beta\rightarrow\kappa$.
For every pair $(i,j)\in\kappa\times\kappa$, we define a $2$-bounded normal coloring $g_{ij}:[\kappa^+- \kappa]^2\rightarrow\kappa^+$ such that if $\alpha<\gamma<\beta$, $f_\beta(\alpha)=i $ and $f_\gamma(\alpha)=j$ then $g_{ij}(\alpha,\beta)=g_{ij}(\gamma,\beta)$. Such $g_{ij}$ can be easily constructed by recursion.

%The definition of $g_{ij}$ is done by induction on $\beta\in\kappa^+$, thus arriving at $\beta$ we define $g_{ij}(\alpha,\beta)$ for every $\alpha\in\beta$.
%We find the unique ordinals $\alpha,\gamma\in\beta$ so that $f_\beta(\alpha)=i$ and $f_\gamma(\alpha)=j$ (the uniqueness follows from the fact that $f_\beta$ is one-to-one) and we choose a new color $c_\beta$.
%Now we let $g_{ij}(\alpha,\beta)=g_{ij}(\gamma,\beta)=c_\beta$.
%Next, we choose $\kappa$ new colors and we define $g_{ij}$ to be one-to-one on all the pairs of the form $(\delta,\beta)$ for every $\delta\notin\{\alpha,\gamma\}$.
%At the end of the inductive process, $g_{ij}$ will be as required.

By our assumption, for every $(i,j)\in\kappa\times\kappa$ there is a club $C_{ij}$ of $\kappa^+$ that is $g_{ij}$-rainbow.
Let $C=\bigcap\{C_{ij}:(i,j)\in\kappa\times\kappa\}$, so $C$ is a club subset of $\kappa^+$ as well.
Pick up any triple $\alpha,\beta,\gamma\in C$ such that $\alpha<\gamma<\beta$.
Set $i=f_\beta(\alpha)$ and $j=f_\gamma(\alpha)$.
It follows that $g_{ij}(\alpha,\beta)=g_{ij}(\gamma,\beta)$, by the definition of $g_{ij}$.
However, $g_{ij}$ is supposed to be one-to-one on $C$ since $C\subseteq C_{ij}$. We have arrived at a contradiction.

\hfill \qedref{clmomega1}

%The above claim works for every successor cardinal, so we need Theorem \ref{thmomega2closed} for regular limit cardinals.
%On the other hand, this theorem requires $\kappa\geq\omega_2$ so we need the above claim for $\omega_1$.
%Together, the statement holds at every uncountable regular cardinal.
\begin{remark}
The key point in the proof of Lemma \ref{clmomega1} is that the club filter over $\kappa^+$ is $\kappa^+$-complete. Any such filter can replace the club filter in the proof.
%and the proof will work by choosing any $\omega_1$-closed filter.
%For example, if we fix a stationary set $T\subseteq\omega_1$ and define $\mathscr{F}_T=\{C\cap T:C$ is a club$\}$, then $\omega_1\nrightarrow^{\rm poly}(\mathscr{F}_T)^2_{\omega_1,2-bd}$.
%This doesn't imply, of course, that $\omega_1\nrightarrow^{\rm poly}(\omega_1-st)^2_{\omega_1,2-bd}$.
\end{remark}

Proposition \ref{constraint} now follows from Lemma \ref{clmomega1} and Lemma \ref{thmomega2closed}.

\begin{proposition}
\label{propomegamm}
If $\omega_2\to^{\rm{poly}} (\omega_1-cl)^2_{<\omega-bdd}$, then 
for every $\sigma\in\omega_1$, $\omega_2\rightarrow^{\rm poly}((\omega_1+\sigma)-cl)^2_{<\omega-bdd}$.
\end{proposition}

\par\noindent\emph{Proof}. \newline 
Fix a countable limit ordinal $\sigma\in\omega_1$ and a $<\omega$-bounded coloring $f:[\omega_2]^2\rightarrow\omega_2$. Apply the hypothesis, we can find two $f$-rainbow closed sets $C_0<C_1$, both of order type $\omega_1$. For each $\beta\geq \sup C_0$, there exists a club $D_\beta \subset C_0$ such that $f\restriction [D_\beta]^2\times \{\beta\}$ is injective. The reason is that for each $\alpha\in C_0$, since the color class $t_{\alpha\beta}=_{def}\{\eta: f(\eta, \beta)=f(\alpha,\beta)\}$ is finite, there exists $\gamma_\alpha\in C_0$ such that $C_0-\gamma_\alpha$ is disjoint from $t_{\alpha\beta}$. Then the elements from $C_0$ that are closure points for $\alpha\mapsto \gamma_\alpha$ will form $D_\beta$. Similarly, for any $\beta_0, \beta_1 \geq \sup C_0$, there is a club $E_{\beta_0, \beta_1}\subset C_0$ disjoint from $t_{\beta_0\beta_1}$.

Let $\delta=\sup C_0$. For each $\alpha\in \lim C_1$, there exists $\beta_\alpha\in C_1\cap \alpha$ such that $t_{\delta\alpha}$ is disjoint from $\alpha-\beta_\alpha$ since $t_{\delta\alpha}$ is finite. Apply Fodor's lemma within $C_1$, we get a stationary $S\subset C_1$ in $\sup C_1$ and some $\xi\in C_1$ such that for all $\alpha\in S$, $\alpha-\xi$ is disjoint from $t_{\delta\alpha}$. As $S$ is essentially a stationary subset of $\omega_1$, it is \emph{fat}. This in particular means we can find a closed subset $F$ of $S$ of order type $\sigma$. Notice that $\{\delta\}\cup F$ is $f$-rainbow. Let $$C^*=\bigcap_{\beta\in F\cup \{\delta\}} D_\beta \cap  \bigcap_{\{\beta_0, \beta_1\} \in [F\cup \{\delta\}]^2} E_{\beta_0,\beta_1},$$ which is a club in $\sup C_0$. It is easy to see that $C^*\cup \{\delta\}\cup F$ is the $f$-rainbow subset as desired.

\hfill \qedref{propomegamm}

In the above theorems, we have seen that prediction principles are very useful when trying to prove negative rainbow statements about \emph{closed} sets. However, parallel negative statements concerning stationary sets cannot be proved in full generality.

\section{Versions of Chang's Conjectures and rainbow sets}\label{CC}

The following boundedness condition naturally arises when versions of Chang's Conjectures are applied to study rainbow Ramsey theory.

\begin{definition}\label{deftbd}
Let $\kappa$ and $\alpha$ be two ordinals.
A function $f:[\kappa]^2\rightarrow\kappa$ is $<\alpha$-\emph{type bounded} iff there is an ordinal $\gamma<\alpha$ so that ${\rm otp}(t_{\alpha\beta})\leq\gamma$ whenever $\alpha<\beta<\kappa$, where $t_{\alpha\beta}=\{\xi\in\beta:f(\xi,\beta)=f(\alpha,\beta)\}$.
\end{definition}

From now on, the domains of the colorings considered are of the form $[C]^2$, where $C$ is a subset of the ordinals.

Type boundedness is not as robust as the cardinality-based boundedness conditions in the sense that if $f$ on $[A]^2$ is $<\alpha$-type bounded, and $h: A\to B$ is a bijection, then $f'$ defined on $[B]^2$ such that $f'(\alpha,\beta)=f(h^{-1}(\alpha), h^{-1}(\beta))$ does not need to be $<\alpha$-type bounded. However, if in addition, we know that $h$ is an order isomorphism, then $f'$ is indeed $<\alpha$-type bounded.

Recall Definition \ref{generalproperty}, $\kappa\rightarrow^{\rm poly}(\theta)^2_{<\alpha-t-bdd}$ iff for every $f:[\kappa]^2\rightarrow\kappa$ which is $<\alpha$-type bounded one can find an $f$-rainbow set of order type $\theta$.
Similar notation like $\kappa\rightarrow^{\rm poly}(\theta-st)^2_{<\alpha-t-bdd}$,  $\kappa\rightarrow^{\rm poly}(\theta-cl)^2_{<\alpha-t-bdd}$ should be interpreted accordingly.
%We shall write $\kappa\rightarrow^{\rm poly}(\theta)^2_{\kappa,\aleph_\eta-bdd}$ if $|t_{\alpha\beta}|\leq\aleph_\eta$ for every $\alpha<\beta<\kappa$.
%The symbol $\kappa\rightarrow^{\rm poly}(\theta)^2_{\kappa,<\gamma-t-bdd}$ means $\kappa\rightarrow^{\rm poly}(\theta)^2_{\kappa,\beta-t-bdd}$ for some fixed $\beta<\gamma$.
%So $\kappa\rightarrow^{\rm poly}(\theta)^2_{\kappa,<\omega_1-t-bdd}$ is weaker than $\kappa\rightarrow^{\rm poly}(\theta)^2_{\kappa,\aleph_0-bdd}$.
%The following statement is essentially due to Todor\v{c}evi\'c:

\begin{definition}\label{partitionrainbow}
Given a boundedness condition $\mathcal{P}$, 
let $$\kappa\to^{\rm{poly}}(par(<\alpha))^2_{\mathcal{P}}$$ abbreviate that for any $f: [\kappa]^2\to \kappa$ satisfying $\mathcal{P}$, there exists a partition $g: \kappa\to \gamma$ for some $\gamma<\alpha$ such that for any $i\in \gamma$, $g^{-1}\{i\}$ is $f$-rainbow.
\end{definition}

\begin{theorem}[Todor\v{c}evi\'c \cite{MR716846}, see also \cite{MR2354904}]\label{pfatypebounded} 
\label{thmtod} 
\begin{enumerate}
\item If \textsf{ZFC} is consistent, so is $\omega_1\to^{\rm{poly}} (par(<\omega_1))^2_{<\omega_1-t-bdd}$.
\item \textsf{PFA} implies $\omega_1\to^{\rm{poly}} (par(<\omega_1))^2_{<\omega_1-t-bdd}$. 
\end{enumerate}
\end{theorem}

\hfill \qedref{thmtod}

\begin{remark}
The theorem proved in \cite{MR716846} was not exactly phrased as the above: it only considers $<\omega$-bounded colorings. However, the proof is essentially the same. The reader is invited to reconstruct the proof either from \cite{MR716846},\cite{MR2354904} or from the proof of Claim \ref{mmproperness} in the later section.
\end{remark}

The following proposition connects Chang's Conjectures with rainbow Ramsey partition relations.

\begin{proposition}
\label{clmcc} Let $\kappa, \lambda, \mu\in \rm{cof}(>\omega), \rho, \theta$ be cardinals. Consider the following statements: 
\begin{enumerate}
\item [$(a)$] $\kappa\to (\mu)_{\lambda, <\rho}^2$, namely for any $f: [\kappa]^2\to \lambda$, there exists $A\subset \kappa$ of order type $\mu$ such that $|f''[A]^2|<\rho$. 
\item [$(b)$] $\kappa\to (\mu-st)_{\lambda, <\rho}^2$, the same as $(a)$ except that we also require $A$ to be stationary in its supremum.
\item [$(c)$] $\mu \to^{\rm{poly}} (\theta)^2_{<\rho-t-bdd}$.
\item [$(d)$] $\mu \to^{\rm{poly}} (par(<{\rm cf}(\mu)))^2_{<\rho-t-bdd}$
\end{enumerate}
Then \begin{enumerate}
\item $(a)+(c)$ implies $\kappa \to^{\rm{poly}} (\theta)^2_{<\lambda^+-t-bdd}$ and
\item $(b)+(d)$ implies $\kappa \to^{\rm{poly}} (\mu-st)^2_{<\lambda^+-t-bdd}$.
\end{enumerate}
\end{proposition}

\par\noindent\emph{Proof}. \newline 
First we show (1).
Assume that $f:[\kappa]^2\rightarrow \kappa$ is $<\lambda^+$-type bounded, as witnessed by $\gamma<\lambda^+$. Recall that we are assuming that $f$ is normal. Let $h: \gamma\to \lambda$ be a bijection.

Define $g:[\kappa]^2\rightarrow \lambda$ by letting $g(\alpha,\beta)=\xi$ iff $\alpha$ is the $h^{-1}(\xi)$-th element of $t_{\alpha\beta}=\{\eta\in\beta:f(\alpha,\beta)=f(\eta,\beta)\}$. Applying $(a)$, one can fix $A\subseteq\kappa$ such that $A$ has order type $\mu$ and $g''[A]^2$ is a subset of $\lambda$ of size $<\rho$, call it $C_A'$. Let $C_A = h^{-1}(C_A')$.

Let $\delta={\rm otp}(C_A)$, so $\delta < \rho$.
We claim that the restriction $f\upharpoonright[A]^2$ is $< \rho$-type bounded as witnessed by $\delta$.
Indeed, given $\alpha<\beta\in A$,  let $s_{\alpha\beta}=t_{\alpha\beta}\cap A$.
Define $h_{\alpha\beta}:s_{\alpha\beta}\rightarrow C_A$ by $h_{\alpha\beta}(\gamma)=h^{-1}(g(\gamma,\beta))$.
It follows from the definition of $g$ that $h_{\alpha\beta}$ is one-to-one and order preserving, hence ${\rm otp}(s_{\alpha\beta})\leq{\rm otp}(C_A)=\delta$ for every $\alpha,\beta\in A$. Therefore, $f\upharpoonright[A]^2$ is $< \rho$-type bounded.
By $(c)$ and the fact that $<\rho$-type-boundedness is preserved under order isomorphisms, we obtain a $f$-rainbow subset $A'$ of $A$ of order type $\theta$.

The proof for (2) is basically the same. By applying $(b)$, one gets $A$ as above satisfying additionally that $A$ is stationary in $\sup A$. Apply $(d)$ instead, we can find a partition of $A$ into $<\rm{cf}(\mu)$ many $f$-rainbow subsets. One of them must have order type $\mu$ and is stationary in its supremum.

\hfill \qedref{clmcc}

\begin{corollary}\label{Changvanilla}
If the Chang's Conjecture and $\omega_1\to^{\rm{poly}}(\omega_1)^2_{<\omega_1-t-bdd}$ hold, then $\omega_2\to^{\rm{poly}}(\omega_1)^2_{\omega-bdd}$.
\end{corollary}
\hfill \qedref{Changvanilla}

\begin{remark}
The consistency of $\omega_2\rightarrow^{\rm poly}(\omega_1)^2_{\omega-bddd}$ does not require large cardinals. In fact, in \cite{MR2354904}, it was shown that it holds in any model obtained by forcing over a model of CH with any poset satisfying the $\omega_1$-covering property. Corollary \ref{Changvanilla} gives a different scenario not covered previously that $\omega_2\rightarrow^{\rm poly}(\omega_1)^2_{\omega-bdd}$ can hold, for example, under Martin's Maximum.
\end{remark}

The following addresses the question raised in \cite{zhang} regarding the rainbow partition relations at the successor of a singular cardinal, as promised in the introduction. 

\begin{corollary}
\label{thmsing} Suppose that:
\begin{enumerate}
\item [$(\aleph)$] $\aleph_{\omega+1}\to (\omega_1-st)_{\aleph_\omega, <\omega_1}^2$ and
\item [$(\beth)$] $\omega_1\rightarrow^{\rm poly}(par (<\omega_1))^2_{<\omega_1-t-bdd}$.
\end{enumerate}
Then $\aleph_{\omega+1}\rightarrow^{\rm poly} (\omega_1-st)^2_{<\aleph_{\omega+1}-t-bdd}$.
\end{corollary}

\hfill \qedref{thmsing}

\begin{remark}
Note that $(\aleph)$ is a consequence of a stationary variation of $(\aleph_{\omega+1},\aleph_\omega)\twoheadrightarrow (\aleph_1, \aleph_0)$.
Here is one way to demonstrate the hypothesis of Theorem \ref{thmsing} is consistent relative to the existence of large cardinals. The argument is essentially due to Levinski, Magidor and Shelah \cite{MR1045371}, while the large cardinal hypothesis needed is improved by Hayut \cite{MR3633795}. By \cite{MR1045371} and \cite{MR3633795}, with the appropriate large cardinal hypothesis (for example, the existence of an $(\omega+1)$-subcompact cardinal), we may suppose in the ground model, GCH holds and there are two strongly inaccessible cardinals $\kappa<\lambda$ satisfying $(\lambda^{+\omega+1}, \lambda^{+\omega})\twoheadrightarrow (\mathrm{stat}(\kappa^{+\omega+1} \cap \mathrm{cof}(\omega)), \kappa^{+\omega})$. More precisely, the symbol means for any regular $\mu\geq \lambda^{+\omega+1}$ and any countable language $L$, there exists $M\prec (H(\mu), L)$ such that $|M\cap \lambda^{+\omega+1}|=\kappa^{+\omega+1}$ and $|M\cap \lambda^{+\omega}|=\kappa^{+\omega}$ and in addition, $M\cap \lambda^{+\omega+1} \cap cof(\omega)$ contains a set $A$ of order type $\kappa^{+\omega+1}$ stationary in its supremum. We argue that the hypothesis of Theorem \ref{thmsing} holds in the forcing extension by $R=_{def}Coll(\omega, \kappa^{+\omega})* P$, where $P$ is a $\lambda$-c.c proper forcing in $V^{Coll(\omega,\kappa^{+\omega})}$ of size $\lambda$ turning $\lambda$ to $\omega_2$ while forcing $\omega_1\to^{\rm poly} (par(<\omega_1))^2_{<\omega_1-t-bdd}$ (see for example the proof of Theorem \ref{pfatypebounded}). Clearly we can assume $R\subset V_\lambda$. Notice that in the forcing extension by $R$, $(\kappa^{+\omega+1})^V = \omega_1$ and $\lambda=\omega_2$.

Let $\dot{c}: [\lambda^{+\omega+1}]^2 \to \lambda^{+\omega}$ be an $R$-name. Let $\mu$ be a sufficiently large regular cardinal and consider the structure $K=(H(\mu), \lambda^{+\omega+1}, \lambda^{+\omega}, R, \dot{c})$. Apply the hypothesis in the ground model, we get $M\prec K$ such that $|M\cap \lambda^{+\omega+1}|=\kappa^{+\omega+1}$ and $|M\cap \lambda^{+\omega}|=\kappa^{+\omega}$ and some $A\subset M\cap \lambda^{+\omega+1}\cap cof(\omega)$ of order type $\kappa^{+\omega+1}$ stationary in its supremum. Let $G\subset R$ be generic over $V$. 

First notice that in $V[G]$, $A$ is a subset of order type $\omega_1$ stationary in its supremum in $V[G]$. The reason is that $Coll(\omega, \kappa^{+\omega})$ satisfies $\kappa^{+\omega+1}$-c.c, which in turn implies that $A$ remains stationary in its supremum whose cofinality becomes $\omega_1$ in $V^{Coll(\omega,\kappa^{+\omega})}$. As $P$ is proper in $V^{Coll(\omega, \kappa^{+\omega})}$, $A$ remains stationary in $\sup A$ in $V[G]$.

Finally, we note that $|M[G]\cap \lambda^{+\omega}|\leq |\kappa^{+\omega}|=\aleph_0$, which in turn implies $c''[A]^2$ is countable since $c'' [A]^2 \subset c''[M[G]\cap \lambda^{+\omega+1}]^2 \subset M[G]\cap \lambda^{+\omega}$. Since $\lambda^{+\omega}$ is singular of countable cofinality, it suffices to show for each $n\in \omega$, $M[G]\cap \lambda^{+n}$ is countable. Since $R$ is of size $\lambda$ satisfying $\lambda$-c.c, a simple calculation reveals that the collection $T(\lambda^{+n})$ of nice $R$-names for an ordinal in $\lambda^{+n}$ has size $\lambda^{\lambda^{+n}} < \lambda^{+\omega}$, by the GCH assumption. Therefore, in $V$, $M\models |T(\lambda^{+n})|<\lambda^{+\omega}$, in particular, $|T(\lambda^{+n})\cap M|\leq \kappa^{+\omega}$. Hence in $V[G]$, we have that $M[G]\cap \lambda^{+n}$ is countable.
\end{remark}

Proposition \ref{clmcc} shows that $\omega_2\rightarrow^{\rm poly}(\omega_1)^2_{\omega-bdd}$ is consistent. A natural question is whether we can enlarge the rainbow set found. The rest of this section is dedicated to the proof of Theorem \ref{mainstat}, which basically says it is consistent to ask for a rainbow subset of order type $\omega_1$ stationary in its supremum. 

\subsection{The role and property of the stationary Chang's Conjecture}
\begin{definition}
\label{defstatcc} 
We shall say that the \emph{stationary Chang's Conjecture} holds iff for every $e:[\omega_2]^{<\omega}\rightarrow\omega_1$ one can find $A\subseteq\omega_2$ such that $\rm{otp}(A)=\omega_1$, $A$ is stationary in $\sup A$ and $e''[A]^{<\omega}$ is countable.
\end{definition}

Our definition is phrased with respect to colorings of all finite subsets, but for all the results in this section we need just the parallel statements applied to colorings of pairs.
It is known that Chang's Conjecture is indestructible under $ccc$ forcing notions.
Let us show that stationary Chang's Conjecture is preserved by $ccc$ forcing notions as well.

\begin{claim}
\label{clmstatcc} Suppose that Stationary Chang's conjecture holds in $V$ and $\mathbb{P}$ is a $ccc$ forcing notion.
If $G\subseteq\mathbb{P}$ is $V$-generic then $V[G]$ satisfies stationary Chang's conjecture.
\end{claim}

\par\noindent\emph{Proof}. \newline 
Let $\name{e}:[\omega_2]^{<\omega}\rightarrow\omega_1$ be a $\mathbb{P}$-name.
For every $t\in[\omega_2]^{<\omega}$ let $\mathcal{A}_t$ be a maximal antichain which determines the possible values of $\name{e}(t)$.
By our assumptions on $\mathbb{P}$, each $\mathcal{A}_t$ is countable.
Define $d:[\omega_2]^{<\omega}\rightarrow\omega_1$ in $V$ by letting $d(t)=\sup\{\alpha\in\omega_1:\exists p\in\mathcal{A}_t,p\Vdash\name{e}(t)=\check{\alpha}\}$.
Since each $\mathcal{A}_t$ is countable we see that $d(t)\in\omega_1$ for every $t$, so $d$ is well-defined.

Apply stationary Chang's Conjecture to $d$ in the ground model, and choose $A\subseteq\omega_2$ of order type $\omega_1$ so that $A$ is stationary in $\sup(A)$ and $d''[A]^{<\omega}$ is countable.
Fix $\delta\in\omega_1$ which bounds $d''[A]^{<\omega}$ and notice that $\Vdash_{\mathbb{P}}\sup\name{e}''[A]^{<\omega}\leq\delta$ as well.
Since $\mathbb{P}$ is $ccc$ we see that $A$ remains stationary in its supremum in $V[G]$, hence the proof is accomplished.

\hfill \qedref{clmstatcc}

\begin{remark}\label{togetStatVersion}
Proposition \ref{clmcc} with $\kappa=\omega_2$ and $\lambda=\mu=\rho=\omega_1$ implies that if the stationary Chang's Conjecture and $\omega_1\to^{\rm{poly}} (par(<\omega_1))^2_{<\omega_1-t-bdd}$ both hold, then $\omega_2\to^{\rm{poly}}(\omega_1-st)^2_{<\omega_2-t-bdd}$ holds, in particular, $\omega_2\to^{\rm{poly}}(\omega_1-st)^2_{\omega-bdd}$ holds.
\end{remark}

For the discussion to follow, it will be convenient to use $$\omega_1\to^{\rm{poly}}_{\aleph_1-ccc} (par(<\omega_1))^2_{<\omega_1-t-bdd}$$ to abbreviate: for any $\aleph_1$-sized $ccc$ poset $Q$ and any $Q$-name $\name{f}$ for a $<\omega_1$-type bounded coloring on $[\omega_1]^2$, there exists a $Q$-name $\name{R}$ for a $ccc$ forcing that adds a partition of $\omega_1$ into countably many $\name{f}$-rainbow sets.

\subsection{Forcing $\omega_1\to^{\rm{poly}}_{\aleph_1-ccc} (par(<\omega_1))^2_{<\omega_1-t-bdd}$ with finite/Easton mixed-support iterations}

We first review the properties of the forcing used for Theorem \ref{pfatypebounded} from \cite{MR716846} and then we modify it for our purpose.

\begin{definition}[Jensen's fast club forcing]
Let $\mathbb{J}_{\omega_1}$ be the poset consisting of $(\sigma, C)$ where $\sigma$ is a countable closed subset of $\omega_1$ and $C$ is a club subset of $\omega_1$. The order is: $(\tau, D)\leq (\sigma, C)$ iff $D\subset C$, $\tau$ end-extends $\sigma$ and $\tau-\sigma\subset C$.
\end{definition}

Let us summarize some properties of this forcing: 
\begin{enumerate}
\item $|\mathbb{J}_{\omega_1}|=2^{\omega_1}$,
\item $\mathbb{J}_{\omega_1}$ is countably closed, 
\item $\mathbb{J}_{\omega_1}$ collapses the continuum of the ground model to $\omega_1$ and
\item $\mathbb{J}_{\omega_1}$ adds a club in $\omega_1$ that is almost contained (except for an initial segment) in any club from the ground model.
\end{enumerate}

\begin{theorem}[Magidor \cite{MR683153}]\label{Magidor}
Let $\kappa$ be regular and $\lambda\geq \kappa$. If $\mathbb{P}$ is a separative forcing poset such that $\mathbb{P}$ is $\kappa$-closed, $|\mathbb{P}|=\lambda$, every condition in $\mathbb{P}$ has $\lambda$ incompatible extensions and $\mathbb{P}$ adds a surjection from $\kappa$ to $\lambda$. Then $\mathbb{P}$ is equivalent to $Coll(\kappa,\lambda)$.
\end{theorem}

As a result, if $V$ satisfies $2^\omega=2^{\omega_1}=\omega_2$, then by Theorem \ref{Magidor}, $\mathbb{J}_{\omega_1}$ is forcing equivalent to $Coll(\omega_1, \omega_2)$.

We will make use of the following theorem, essentially due to Todor\v{c}evi\'{c} \cite[Lemma 1]{MR716846}.

\begin{lemma}
\label{lemtod1} Assume that:
\begin{enumerate}
\item [$(a)$] $\mathbb{Q}$ is a $ccc$ forcing notion in $V$.
\item [$(b)$] $\name{f}:[\omega_1]^2\rightarrow\omega_1$ is a $\mathbb{Q}$-name for a $<\omega_1$-type bounded coloring.
\item [$(c)$] The set $(G_0\ast G_1)\times G_{\mathbb{Q}}$ is $V$-generic with respect to the forcing notion $(Add(\omega,\omega_1)\ast\mathbb{J}_{\omega_1})\times\mathbb{Q}$.
\end{enumerate}
Then in $V[(G_0\ast G_1)\times G_{\mathbb{Q}}]$, there is a $ccc$ forcing notion which adds a partition of $\omega_1$ into countably many $f$-rainbow subsets.
\end{lemma}

\hfill \qedref{lemtod1}

Let $E$ be the class of even ordinals and $O$ be the class of odd ordinals. Given $p\in Coll(\omega_1, \omega_2)$, let the \emph{height} of $p$ be denoted as $ht(p)=\sup \mathrm{dom}(p)$.

The main idea behind the proof of the following theorem is to use a mixed support iteration of Cohen forcing and Lévy Collapse.
This was the strategy of Mitchell in his classical proof of the consistency of the tree property at $\aleph_2$, see \cite{MR0313057}. 

The support for the forcing iterations used in \cite{MR716846} and \cite{MR0313057} is finite on the Cohen part and countable on the countably closed part. We will use \emph{Easton support} on the countably closed part and in addition, we require the conditions to be of certain shape. This will be useful for the lifting argument in Subsection \ref{lifting}.

Here is a short history of how these ideas originated: 
Silver, in his proof of the consistency of the Chang's Conjecture from a $\omega_1$-Erd\H{o}s cardinal, invented the \emph{Silver Collapse}, which is a modification of the usual Lévy Collapse requiring each condition to be of a specific shape. The method was later adapted by Kunen \cite{MR495118} to establish the consistency of the existence of a saturated ideal ideal on $\omega_1$ from the existence of a huge cardinal.
Then Laver \cite{MR673792} modified Kunen's method further, using 
the Easton support, to establish the consistency of an  $(\aleph_2, \aleph_2, \aleph_0)$-saturated ideal on $\omega_1$.

A set $A\subseteq\lambda$ is an \emph{Easton set} iff $|A\cap\rho|<\rho$ for every limit regular cardinal $\rho<\lambda$. The following theorem was basically proved in \cite{MR716846}. We include the proof adapted for the modified context for completeness.

\begin{lemma}
\label{thmmahlo} If $\lambda$ is $2$-Mahlo then there is a forcing notion $\mathbb{S}$ with the following properties:
\begin{enumerate}
\item [$(a)$] $\mathbb{S}$ is proper and $\lambda$-cc.
\item [$(b)$] $\mathbb{S}$ forces $\lambda=\omega_2$.
\item [$(c)$] $\mathbb{S}$ forces $\omega_1\to^{\rm{poly}}_{\aleph_1-ccc} (par(<\omega_1))^2_{<\omega_1-t-bdd}$.
\end{enumerate}
\end{lemma}

\par\noindent\emph{Proof}. \newline 
Let $E$ be the class of even ordinals and let $O$ be the class of odd ordinals.
We define a mixed support iteration $(\mathbb{S}_\alpha,\name{\mathbb{T}}_\beta: \alpha\leq\lambda,\beta<\lambda)$ as follows:
\begin{enumerate}
\item [$(\alpha)$] If $\alpha\in E$ then $\Vdash_{\mathbb{S}_\alpha}\name{\mathbb{T}}_\alpha=Add(\omega,\omega_1)$.
\item [$(\beta)$] If $\alpha\in O$ then $\Vdash_{\mathbb{S}_\alpha}\name{\mathbb{T}}_\alpha=Coll(\omega_1,\omega_2)$.
\item [$(\gamma)$] If $\alpha$ is a limit regular cardinal, then $\mathbb{S}_\alpha$ is the direct limit of $(\mathbb{S}_\beta:\beta\in\alpha)$.
\item [$(\delta)$] If $\alpha$ is not a limit regular cardinal, then $t\in\mathbb{S}_\alpha$ iff $t\upharpoonright\beta\in\mathbb{S}_\beta$ for every $\beta\in\alpha, |\{\gamma\in E\cap\alpha:t(\gamma)\neq\varnothing\}|<\aleph_0$ and there exists $\nu\in\omega_1$ such that for every $\gamma\in O\cap\alpha, \Vdash_{\mathbb{S}_\gamma}ht(t(\gamma))<\nu$.
\end{enumerate}

If $t\in\mathbb{S}_\alpha$, in the case of $(\delta)$, then the least $\nu\in\omega_1$ which bounds $ht(t(\gamma))$ for every $\gamma\in O\cap\alpha$ will be denoted by $ht(t)$.
We shall say that $t$ has \emph{bounded height}.
Property $(\delta)$ implies that if $\alpha\leq\lambda, t$ belongs to the inverse limit of $(\mathbb{S}_\beta:\beta\in\alpha)$ and $t$ is of bounded height then $t\in\mathbb{S}_\alpha$ iff $O\cap{\rm supp}(t)$ is an Easton set and $E\cap{\rm supp}(t)$ is finite.

We claim that $\mathbb{S}_\lambda=\mathbb{S}$ satisfies all the requirements of the theorem. By design, it is easy to see $\mathbb{S}$ satisfies $\lambda$-cc and forces $\lambda=\omega_2$.
Let us show that $\mathbb{S}_\alpha$ is proper for every $\alpha\leq\lambda$.

Firstly, for every $\alpha\leq\lambda$ let $D_\alpha=\{p\in\mathbb{S}_\alpha: \forall\gamma\in E\cap\alpha,\exists s_\gamma\in Add(\omega,\omega_1), p\upharpoonright\gamma\Vdash_{\mathbb{S}_\alpha} \check{s}_\gamma=p(\gamma)\}$.
By induction on $\alpha$ it follows that $D_\alpha$ is dense (in $\mathbb{S}_\alpha$) for every $\alpha\leq\lambda$.
Hence we may assume, without loss of generality, that all our conditions belong to the pertinent $D_\alpha$.
For every $\alpha\leq\lambda$, one can see that $\mathbb{S}_\alpha$ projects onto $Add(\omega,\omega_1\times(E\cap\alpha))$.
Each element $q$ of $Add(\omega,\omega_1\times(E\cap\alpha))$ is a function from a finite subset of $E\cap\alpha$ to $Add(\omega,\omega_1)$, ordered by reverse inclusion.
Hence for $\mathbb{S}$ can be identified with the two-step iteration $Add(\omega,\omega_1\times (\lambda\cap E))\ast \mathbb{S}/Add(\omega,\omega_1\times(\lambda\cap E))$.

Fix $\alpha\leq\lambda$.
Given $p\in\mathbb{S}_\alpha$ we denote the Cohen part of $p$ by $\sigma(p)$.
Explicitly, the domain of $\sigma(p)$ is ${\rm supp}(p)\cap E$ and if $\gamma\in{\rm supp}(p)\cap E$ then $\sigma(p)(\gamma)\in Add(\omega,\omega_1)$ and $p\upharpoonright\gamma\Vdash_{\mathbb{S}_\gamma} \sigma(p)(\gamma)=p(\gamma)$.
The main disadvantage of the Cohen part is the lack of completeness, but we can define the direct extension $\leq_\alpha^*$ on $\mathbb{S}_\alpha$ by $q\leq_\alpha^*p$ iff $q\leq_{\mathbb{S}_\alpha}p$ and $\sigma(q)=\sigma(p)$.
Observe that $\leq_\alpha^*$ is $\omega_1$-closed. In fact, any countable decreasing $\leq_\alpha^*$-sequence has a greatest lower bound.

If $p\in\mathbb{S}_\alpha$ and $t \leq_{Add(\omega,\omega_1\times(E\cap \alpha))}\sigma(p)$ then we define an amalgamation $p\wedge t$ of $p$ and $t$ as follows.
Let $q=p\wedge t$ where ${\rm dom}(q)$ is ${\rm dom}(p)\cup{\rm dom}(t)$ and if $\gamma\in \rm{dom}(p)-\rm{dom}(t)$, then $q(\gamma)=p(\gamma)$ and if $\gamma\in \rm{dom}(t)$, $q(\gamma)=t(\gamma)$.
Notice that $q\leq_{\mathbb{S}_\alpha}p,t$.
We proceed now to properness.
The following comes from \cite{MR716846}, but we unfold the proof for completeness.

\begin{claim}
\label{clmtod} Suppose that $\alpha\leq\lambda, p\in\mathbb{S}_\alpha$ and $\name{\beta}$ is an $\mathbb{S}_\alpha$-name of an ordinal.
Then there exists $q\in\mathbb{S}_\alpha$ such that $q\leq_\alpha^*p$ and there are a maximal antichain $\{t_i:i\in\omega\}\subseteq Add(\omega,\omega_1\times(\alpha\cap E))$ below $\sigma(p)$ and $B=\{\beta_i:i\in\omega\}\subseteq{\rm Ord}$ such that $q\wedge t_i\Vdash_{\mathbb{S}_\alpha}\name{\beta}=\beta_i$ for every $i\in\omega$.
\end{claim}

\par\noindent\emph{Proof}. \newline 
By induction on $i\in \omega_1$, we define $p_i\in\mathbb{S}_\alpha$,  $t_i\in Add(\omega,\omega_1\times\alpha)$ and find some $\delta<\omega_1$ such that $\{t_i:i\in\delta\}$ is a maximal antichain below $\sigma(p)$ and if $i<j<\delta$ then $p_j\leq_\alpha^*p_i$.

Arriving at $\gamma<\omega_1$, if $\{t_i:i\in\gamma\}$ is a maximal antichain below $\sigma(p)$ then we let $\delta=\gamma$ and we are done.
If not, let $q$ be the greatest lower bound of $(p_i:i\in\gamma)$ with respect to $\leq_\alpha^*$, and let $t_\gamma \leq \sigma(p)$ be such that $t_i\perp t_\gamma$ for every $i\in\gamma$.
Choose $r\in\mathbb{S}_\alpha$ such that $r \leq_{\mathbb{S}_\alpha} q\wedge t_\gamma$ and $r$ forces a value to $\name{\beta}$, say $\beta_\gamma$.
From $r$ we can extract $p_\gamma,s_\gamma$ such that $p_{\gamma}\leq_\alpha^* q$ and $s_\gamma\leq_{Add(\omega,\omega_1\times(E\cap \alpha))}t_\gamma$.
Notice that $p_\gamma\wedge s_\gamma\Vdash\name{\beta}=\beta_\gamma$.
Since $Add(\omega,\omega_1\times(E\cap \alpha))$ is $ccc$, the process must be terminated at some $\delta\in\omega_1$, so letting $B=\{\beta_i:i\in\delta\}$ we are done.

\hfill \qedref{clmtod}

We conclude from the above claim that each $\mathbb{S}_\alpha$ is proper.
In particular, $\mathbb{S}$ is proper, thus $(a)$ and $(b)$ of our theorem are established and we can proceed to $(c)$.
\begin{claim}[Page 713-714, \cite{MR716846}]
\label{clmpreserveccc} 
If $\mathbb{Q}\in V$ is $ccc$ then $\Vdash_{\mathbb{S}}\mathbb{Q}$ is $ccc$.
\end{claim}
\hfill \qedref{clmpreserveccc}

Observe that if $\alpha<\lambda$ is Mahlo then $\mathbb{S}_\alpha$ forces that $\mathbb{S}_{[\alpha,\lambda)}$ is forcing equivalent to $\mathbb{S}$ as defined in $V^{\mathbb{S}_\alpha}$. An easy modification of Theorem 5.4 in Baumgartner's survey on iterated forcing \cite{MR823775} will give a proof.

To establish $(c)$, suppose that $\mathbb{Q}$ is a $ccc$ forcing notion of size $\aleph_1$ and $\name{f}:[\omega_1]^2\rightarrow\omega_1$ is a $\mathbb{Q}$-name of a $<\omega_1$-type bounded function.
Let $H\subseteq\mathbb{S}$ be $V$-generic.
Find a sufficiently large Mahlo cardinal $\alpha<\lambda$ such that $\mathbb{Q},\name{f}\in V[H\upharpoonright\alpha]$.

Consider the two-step iteration after the $\alpha$th stage.
By Theorem \ref{Magidor}, we see that in $V[H\upharpoonright\alpha]$ these two steps are equivalent to $Add(\omega,\omega_1)\ast\mathbb{J}_{\omega_1}$.
Applying Lemma \ref{lemtod1} we can find in $V[H\upharpoonright\alpha+2]$ a $\mathbb{Q}$-name $\name{\mathbb{R}}$ of a $ccc$ forcing notion which adds a countable partition of $\omega_1$ into $\name{f}$-rainbow sets.
Since $\mathbb{S}_{[\alpha+2,\lambda)}$ is equivalent to $\mathbb{S}$ we conclude by Claim \ref{clmpreserveccc} that $\mathbb{Q}\ast\name{\mathbb{R}}$ is $ccc$ in $V[H]$, thus proving $(c)$.

\hfill \qedref{thmmahlo}

We conclude this subsection by analyzing the term forcing associated with $\mathbb{S}$. Let $\mathbb{R}_\lambda$ be a forcing notion consisting of functions on $\lambda\cap O$. For each $f\in \mathbb{R}_\lambda$, the following are satisfied: 
\begin{enumerate}
\item for any $\alpha\in \lambda\cap O$, $\Vdash_{\mathbb{S}_\alpha} f(\alpha)\in Coll(\omega_1,\omega_2)$;
\item the support of a condition $f\in\mathbb{R}_\lambda$, namely $\rm{supp}(f)=_{def}\{\alpha\in \lambda\cap O: f(\alpha)\neq \emptyset\}$, is an Easton set;
\item there exists $\delta\in\omega_1$ such that $\Vdash_{\mathbb{S}_\alpha}ht(f(\alpha))<\delta$ for every $\alpha\in{\rm supp}(f)$. We denote by $ht(f)$ the minimal $\delta\in\omega_1$ which satisfies this property.
\end{enumerate}

If $f,g\in\mathbb{R}_\lambda$ then $g\leq_{\mathbb{R}_\lambda}f$ iff ${\rm supp}(f)\subseteq{\rm supp}(g)$ and if $\alpha\in{\rm supp}(f)$ then $\Vdash_{\mathbb{S}_\alpha}f(\alpha)\leq g(\alpha)$.

By the nature of the collapse and the maximality principle of forcing, if $f,g\in\mathbb{R}_\lambda$ are compatible then there is a greatest lower bound $h\in\mathbb{R}_\lambda$ for $f,g$. 

\begin{claim}
\label{clmproj} Let $\mathbb{D}_\lambda= Add(\omega,\omega_1\times(\lambda\cap E))\times\mathbb{R}_\lambda$. \newline 
Then $\mathbb{D}_\lambda$ projects onto $\mathbb{S}$.
\end{claim}

\par\noindent\emph{Proof}. \newline 
Suppose that $(s,f)\in\mathbb{D}_\lambda$.
Define $p=p_{s,f}\in\mathbb{S}$ as follows.
Let ${\rm supp}(p)={\rm supp}(s)\cup{\rm supp}(f)$.
If $\gamma\in{\rm supp}(s)$ let $p(\gamma)=s(\gamma)$ and if $\gamma\in{\rm supp}(f)$ let $p(\gamma)=f(\gamma)$.
Define $\pi:\mathbb{D}_\lambda\rightarrow\mathbb{S}$ by $\pi((s,f))=p_{s,f}$.
We claim that $\pi$ projects $\mathbb{D}_\lambda$ to $\mathbb{S}$.

Clearly, $\pi$ sends $0_{\mathbb{D}_\lambda}$ to $0_{\mathbb{S}}$ and preserves order.
Suppose that $(s,f)\in\mathbb{D}_\lambda$ and $q\in\mathbb{S}$ satisfies $p_{s,f}\geq_{\mathbb{S}}q$.
We try to find $(t,g)\in\mathbb{D}_\lambda$ such that $(s,f)\geq_{\mathbb{D}_\lambda}(t,g)$ and $q\geq_{\mathbb{S}}\pi(t,g)=p_{t,g}$.
Let $t=\sigma(q)$.
Define ${\rm supp}(g)=O\cap{\rm supp}(q)$.
If $\alpha\in{\rm supp}(g)$ then let $g(\alpha)$ be an $\mathbb{S}_\alpha$-term for which if $q\upharpoonright\alpha$ belongs to the generic subset of $\mathbb{S}_\alpha$, then $g(\alpha)=q(\alpha)$ and if not then $g(\alpha)=f(\alpha)$.

One can verify that $g$ is a well-defined condition, in particular if $\nu=\max\{ht(f),ht(q)\}\in\omega_1$ then $ht(g)\leq\nu$ and hence $g\in\mathbb{R}_\lambda$.
Since $q\upharpoonright\alpha\Vdash_{\mathbb{S}_\alpha}f(\alpha)\geq q(\alpha)$ for every $\alpha\in{\rm supp}(g)$, we have $f\geq_{\mathbb{R}_\lambda}g$.
To see that $q\geq_{\mathbb{S}}p_{t,g}$, induct on $\alpha\geq\lambda$.
Arriving at $\alpha$, assume that $q\upharpoonright\alpha\geq_{\mathbb{S}_\alpha}p_{t,g}\upharpoonright\alpha$.
Now if $\alpha\in E$, use the fact that $t=\sigma(q)$ for the $\alpha$-th stage.
If $\alpha\in O$, by the induction hypothesis, $p_{t,g} \leq q\restriction \alpha$ and the fact that $p_{t,g}(\alpha)=g(\alpha)$, we know that $q\upharpoonright\alpha\Vdash_{\mathbb{S}_\alpha}g(\alpha)=q(\alpha)$.
In all other cases use the definition of the iteration.

\hfill \qedref{clmproj}

Clearly $Add(\omega, \omega_1\times \lambda)$ is isomorphic to $Add(\omega, \omega_1 \times (\lambda\cap E))$. In what follows, they are identified to be the same.
Notice that $\mathbb{S}$ and $\mathbb{R}$ have uniform definitions with the parameter $\lambda$ and the $\omega_1$ in the ground model. 
Therefore, in any model of set theory with $\nu$ being a 2-Mahlo cardinal, we can make sense of $\mathbb{S}_\nu$ and $\mathbb{R}_\nu$ respectively.

\subsection{Getting the stationary Chang's Conjecture}\label{lifting}

Let $pr:[{\rm Ord}]^2\leftrightarrow{\rm Ord}$ be a pairing function which satisfies the requirement that if $\delta=pr(\beta,\gamma)$ then $\beta\leq \gamma\leq\delta$. Recall that a cardinal $\kappa$ is \emph{huge} if there exists an elementary embedding $\jmath: V\to M$ with critical point $\kappa$ and ${}^{\jmath (\kappa)}M \subset M$. 

\begin{lemma}
\label{thmstatcc} Assume that:
\begin{enumerate}
\item [$(a)$] $\kappa$ is a huge cardinal, $\jmath$ is a huge embedding for $\kappa$ and $\lambda=\jmath(\kappa)$.
\item [$(b)$] $\mathbb{P} = \langle \mathbb{P}_\alpha, \name{\mathbb{Q}}_\beta: \alpha\leq \kappa, \beta<\kappa \rangle$ is a finite support iteration of length $\kappa$.
\item [$(c)$] $\mathbb{P}_0=Coll(\omega,<\kappa)$.
\item [$(d)$] If $\alpha\in(0,\kappa), pr^{-1}(\alpha)=(\alpha_0,\alpha_1)$ and $\mathbb{P}_{\alpha_1}\cap V_{\alpha_0}$ is a regular suborder of $\mathbb{P}_{\alpha_1}$ then $\mathbb{P}_{\alpha+1}=\mathbb{P}_\alpha\ast\name{\mathbb{Q}}_\alpha$ where $\name{\mathbb{Q}}_\alpha$ is a $\mathbb{P}_\alpha$-name of $(\mathbb{R}_\kappa)^{V^{\mathbb{P}_{\alpha_1}\cap V_{\alpha_0}}}$.
\item [$(e)$] In all other cases $\name{\mathbb{Q}}_\alpha$ is a name of the trivial forcing.
\end{enumerate}
Then $\mathbb{P}$ satisfies $\kappa$-cc and $V^{\mathbb{P}\ast\mathbb{S}}$ satisfies the stationary Chang's Conjecture, where $\mathbb{S}$ is the $\mathbb{P}$-name for $\mathbb{S}_\lambda$, as defined in Lemma \ref{thmmahlo}.
\end{lemma}

\par\noindent\emph{Proof}. \newline 
To see $\mathbb{P}$ satisfies $\kappa$-cc, work by induction on $\alpha\leq\kappa$. Assume that $(\alpha_0,\alpha_1)$ satisfy $\alpha=pr(\alpha_0,\alpha_1)$.
By the induction hypothesis we know that $\mathbb{P}_\alpha$ is $\kappa$-cc.
In the trivial stages it is clear that $\mathbb{P}_\alpha$ forces that $\name{\mathbb{Q}}_\alpha$ is $\kappa$-cc, so assume that we arrived at a non-trivial stage, namely $\mathbb{P}_{\alpha_1}\cap V_{\alpha_0}$ is a regular suborder of $\mathbb{P}_{\alpha_1}$ and hence of $\mathbb{P}_\alpha$.
From the fact that $\mathbb{P}_\alpha$ is $\kappa$-cc we infer that $\Vdash_{\mathbb{P}_{\alpha_1}\cap V_{\alpha_0}} \mathbb{P}_\alpha/\mathbb{P}_{\alpha_1}\cap V_{\alpha_0}$ is $\kappa$-cc, and that $\Vdash_{\mathbb{P}_{\alpha_1}\cap V_{\alpha_0}}(\mathbb{R}_\kappa)^{V^{\mathbb{P}_{\alpha_1}\cap V_{\alpha_0}}}$ is $\kappa$-cc.
As $\mathbb{P}_{\alpha_1}\cap V_{\alpha_0}$ is a small forcing, it preserves the weak compactness of $\kappa$.
Hence in $V^{\mathbb{P}_{\alpha_1}\cap V_{\alpha_0}}$ we see that $\kappa$-cc forcings are $\kappa$-Knaster and we are done.

Next we show that $V^{\mathbb{P}\ast\mathbb{S}}$ is a model of the stationary Chang's Conjecture.
Observe that by the construction, $\mathbb{P}\ast\mathbb{R}_\lambda$ is a regular suborder of $\jmath(\mathbb{P})$.
Let $G\ast H\subseteq\mathbb{P}\ast\mathbb{R}_\lambda$ be $V$-generic, and let $\tilde{G}\subseteq\jmath(\mathbb{P})/G\ast H$ be generic over $V[G\ast H]$.
Apply Silver's criterion to lift $\jmath$ and to get $\jmath_0:V[G]\rightarrow M[\tilde{G}]$.
Let $h\subseteq Add(\omega,\kappa\times \lambda)$ be generic over $V[\tilde{G}]$ and let $L\subseteq\mathbb{S}$ be generic over $V[G]$, as  given by $h\times H$ and Claim \ref{clmproj}. All of the above happens in $V[\tilde{G}][h]$.

Note that in $V[\tilde{G}][h]$, $\jmath_0\restriction Add(\omega, \kappa\times \lambda): Add(\omega, \kappa\times \lambda) \to Add(\omega, \lambda \times \jmath(\lambda))$ is a complete embedding, with the isomorphic image being $Add(\omega, \kappa\times \jmath''\lambda)$. As a result, $\jmath_0 '' h \subset  Add(\omega, \kappa\times \jmath''\lambda)$ is generic over $V[\tilde{G}]$.

We shall argue that by going to a further forcing extension, there is some $\tilde{L}\subseteq\jmath_0(\mathbb{S})$ generic over $M[\tilde{G}]$, so that $\jmath_0''L\subseteq\tilde{L}$.
To this end, notice that $\jmath_0''H\in M[\tilde{G}]$.
Now in $M[\tilde{G}]$, we find a \emph{master condition} for $\jmath_0'' H$. More precisely, we define a lower bound $r\in\jmath_0(\mathbb{R}_\lambda)$ for all the conditions in $\jmath_0''H$.
Let $${\rm supp}(r)=\bigcup\{\jmath_0({\rm supp}(q)):q\in H\}.$$
Each set of the form $\jmath_0({\rm supp}(q))$ is Easton in $M[\tilde{G}]$ by the elementarity of $\jmath_0$ and hence ${\rm supp}(r)$ is also an Easton set in $M[\tilde{G}]$.
Indeed, $M[\tilde{G}]\models\lambda=\omega_1$, so ${\rm supp}(r)$ is a union of $\omega_1$ many Easton sets in $M[\tilde{G}]$ and hence Easton.

For every $\alpha\in{\rm supp}(r)$ let $r(\alpha)$ be a $\jmath_0(\mathbb{S})\upharpoonright\alpha$-term for the greatest lower bound for $\{\jmath_0(q)(\alpha):q\in H\}$.
To verify that $r$ is a condition in $\jmath_0(\mathbb{R}_\lambda)$, let $\gamma_q=ht(q)$ for each $q\in H$, so $\gamma_q\in\kappa$.
Since $\kappa={\rm crit}(\jmath_0)$ we see that $ht(\jmath_0(q))=\gamma_q$ as well.
Hence $ht(r)\leq\kappa$. As $M[\tilde{G}]\models ``\kappa<\lambda = \omega_1$'', we see that $r$ is a condition in $\jmath_0(\mathbb{R}_\lambda)$ indeed.

Let $\tilde{h}\times\tilde{H}\subseteq \jmath_0(Add(\omega,\kappa\times\lambda))/\jmath_0''h\times \jmath_0(\mathbb{R}_\lambda)/r$ be generic over $V[\tilde{G}\ast h]$. In $V[\tilde{G}*h][\tilde{h}\times \tilde{H}]$, 
one can lift $\jmath_0$ to $\jmath_1:V[G][h\times H]\rightarrow M[\tilde{G}][\tilde{h}\times\tilde{H}]$.
Let $\tilde{L}\subseteq\jmath_0(\mathbb{S})$ be $M[\tilde{G}]$-generic as given by $\tilde{h}\times\tilde{H}$ and Claim \ref{clmproj}.
Let $\imath=\jmath_1\upharpoonright V[G\ast L]$. It is clear that $\imath:V[G\ast L]\rightarrow M[\tilde{G}\ast\tilde{L}]$ is an elementary embedding as desired. This elementary embedding lives in $V[\tilde{G}*h][\tilde{h}\times \tilde{H}]$.

We claim that $V[G\ast L]$ satisfies stationary Chang's conjecture.
In order to prove this, suppose that $e:[\lambda]^{<\omega}\rightarrow\kappa, e\in V[G\ast L]$ and recall that $\lambda=\omega_2,\kappa=\omega_1$ in $V[G\ast L]$.
Our goal is to find $A\subseteq\lambda, {\rm otp}(A)=\omega_1$, $A$ is stationary in its supremum such that $e'' [A]^{<\omega}$ is countable.

In $V$, let $A=\jmath''(\lambda\cap{\rm cof}^V(\omega))$, namely the pointwise image of $\jmath$ over the set of $\alpha\in\lambda$ with countable cofinality in $V$.
Let $\beta=\sup(A)$.
Since $\jmath:V\rightarrow M$ is a huge embedding with target $\lambda$, we see that $M$ is $\lambda$-closed so $A\in M$.
Moreover, as $\jmath$ is continuous at points of countable cofinality, $A$ is stationary in $\beta$.

Since $\jmath(\mathbb{P})$ satisfies $\lambda$-cc, we conclude that $A$ is stationary in $\beta$ in $M[\tilde{G}]$.
Recall that $\mathbb{S}$ is proper, so by elementarity $\jmath_0(\mathbb{S})$ is proper in $M[\tilde{G}]$.
Recall that $\lambda=\omega_1$ in $M[\tilde{G}]$, so $A$ remains stationary in $\beta$ in $M[\tilde{G}\ast\tilde{L}]$.

Now $\kappa<\lambda=\omega_1$ in $M[\tilde{G}\ast\tilde{L}]$, and $\imath(e)''[A]^{<\omega}\subseteq\kappa$, so $M[\tilde{G}\ast\tilde{L}]$ knows that there exists a subset of $\imath(\lambda)$ of order type $\omega_1$, stationary in its supremum, with countable image under $\imath(e)$.
By elementarity, this statement is true also in $V[G\ast L]$, as required.

\hfill \qedref{thmstatcc}

Now we can finish the proof of Theorem \ref{mainstat}.

\par\noindent\emph{Proof of Theorem \ref{mainstat}}. \newline 
First, we force with $\mathbb{P}\ast\mathbb{S}$ as in Lemma \ref{thmstatcc} where $\kappa$ is huge witnessed by $\jmath: V\to M$, $\lambda=\jmath(\kappa)$ and $\mathbb{S}=(\mathbb{S}_\lambda)^{V^{\mathbb{P}}}$.
Let $G\subseteq\mathbb{P}\ast\mathbb{S}$ be $V$-generic. Notice that in $V[G]$, $\kappa=\omega_1$ and $2^\omega=2^{\omega_1}=\omega_2=\lambda$.
From Lemma \ref{thmstatcc} we infer that $V[G]$ satisfies the stationary Chang's Conjecture.

Working in $V[G]$, we define a finite support iteration $\mathbb{Q}$ of $ccc$ forcings of length $\omega_2$, forcing $\omega_1\rightarrow^{\rm poly}(par(<\omega_1))^2_{<\omega_1-t-bdd}$. This is possible, since by Lemma \ref{thmmahlo}(c) $\omega_1\rightarrow^{\rm poly}_{\aleph_1-ccc}(par(<\omega_1))^2_{<\omega_1-t-bdd}$ holds in $V[G]$. Using the standard book-keeping trick, one can define an $\omega_2$-length finite support iteration of $\aleph_1$-sized $ccc$ posets to force $\omega_1\rightarrow^{\rm poly}(par(<\omega_1))^2_{<\omega_1-t-bdd}$.

Let $H\subseteq\mathbb{Q}$ be $V[G]$-generic.
Since $\mathbb{Q}$ is $ccc$, we deduce from Claim \ref{clmstatcc} that $V[G][H]$ satisfies the stationary Chang's Conjecture. By Remark \ref{togetStatVersion}, $\omega_2\rightarrow^{\rm poly}(\omega_1-st)^2_{\omega-bdd}$ (even $\omega_2\rightarrow^{\rm poly}(\omega_1-st)^2_{<\omega_2-t-bdd}$) holds in $V[G][H]$.

\hfill \qedref{mainstat}

\section{\textsf{PFA} does not imply $\omega_2\to^{\rm poly} (\omega_1-st)^2_{\omega-bdd}$}\label{pfanotimplies}

Recall Theorem \ref{pfatypebounded} that PFA implies $\omega_1\to^{\rm poly} (\omega_1-st)^2_{<\omega_1-t-bdd}$, hence in particular, $\omega_2\to^{poly} (\omega_1-st)^2_{<\omega_1-t-bdd}$. We have seen that by Theorem \ref{mainstat} that $\omega_2\to^{poly} (\omega_1-st)^2_{\omega-bdd}$ is consistent relative to the existence of a huge cardinal. A natural question is whether \textsf{PFA} already implies $\omega_2\to^{poly} (\omega_1-st)^2_{\omega-bdd}$. In this section, we show that the answer is no.

Let $(\ast)$ be the statement which says that there exists a sequence of sets $(S_\alpha\in[\alpha]^{\leq\aleph_0}:\alpha\in S^{\omega_2}_\omega)$ such that for every $\beta\in S^{\omega_2}_{\omega_1}$ there is a club $c_\beta\subseteq\beta \cap \mathrm{cof}(\omega)$ so that if $\alpha\in c_\beta$ then $c_\beta\cap\alpha\subseteq S_\alpha$.

\begin{claim}
\label{clmstar} If $(\ast)$ holds then $\omega_2\nrightarrow^{\rm poly}(\omega_1-st)^2_{\omega-bdd}$.
\end{claim}

\par\noindent\emph{Proof}. \newline 
Given the sequence $(S_\alpha:\alpha\in S^{\omega_2}_\omega)$ we define a normal $\omega$-bounded coloring $f:[\omega_2]^2\rightarrow\omega_2$ such that for any $\alpha\in S^{\omega_2}_\omega$, $f(\cdot, \alpha)\restriction S_\alpha$ is constant. This is possible since each $S_\alpha$ is countable.

Suppose we are given $A\subseteq\omega_2$ such that ${\rm otp}(A)=\omega_1, \beta=\sup(A)$ and $A$ is stationary in $\beta$.
Choose $\alpha\in A\cap c_\beta \cap S^{\omega_2}_\omega$ such that $|A\cap c_\beta\cap \alpha|\geq 2$. By $(\ast)$, $A\cap c_\beta\cap \alpha \subseteq c_\beta\cap \alpha\subseteq S_\alpha$, so $A$ is not $f$-rainbow.

\hfill \qedref{clmstar}

\begin{remark}
$(*)$ is a consequence of $\square_{\omega_1, \omega}$. Therefore, unlike $\omega_2\rightarrow^{\rm poly}(\omega_1)^2_{\omega-bdd}$, the stationary variation $\omega_2\rightarrow^{\rm poly}(\omega_1-st)^2_{\omega-bdd}$ has non-trivial large cardinal strength.
\end{remark}

Our objective is to show that $(\ast)$ is consistent with \textsf{PFA}.
We begin with a ground model which satisfies \textsf{PFA}, we force $(\ast)$ and we would like to argue that \textsf{PFA} is preserved in the generic extension.
The proof follows Beaudoin \cite{MR1099782} and Magidor, and we refer to \cite[Section 7]{MR2354904} for a similar argument and a comprehensive discussion concerning the essential idea of \cite{MR1099782}.

\begin{theorem}
\label{thmpfa} \textsf{PFA} is consistent with $(\ast)$, and hence does not imply the positive relation $\omega_2\rightarrow^{\rm poly}(\omega_1-st)^2_{\omega-bdd}$.
\end{theorem}

\par\noindent\emph{Proof}. \newline 
We define a forcing notion $\mathbb{P}$ to add witness for $(\ast)$ with initial segments.
A condition $p\in\mathbb{P}$ iff
\begin{itemize}
\item $\mathrm{dom}(p)$ is a proper initial segment of $S^{\omega_2}_\omega$, 
\item if $\alpha\in{\rm dom}(p)$, then $p(\alpha)$ is a countable subset of $\alpha$,
\item for every $\beta\in S^{\omega_2}_{\omega_1}$ and $\beta\leq \sup \mathrm{dom}(p)$, there exists a club $c_\beta\subset \beta$ such that for any $\alpha\in c_\beta$, $c_\beta \cap \alpha\subset g(\alpha)$.
\end{itemize}
 The order $\leq_\mathbb{P}$ is end extension, namely, $q\leq_{\mathbb{P}} p$ iff $\mathrm{dom}(q)\supseteq  \mathrm{dom}(p)$ and $q\restriction \mathrm{dom}(p)=p$.

Observe that $\mathbb{P}$ is countably closed and $(\omega_1+1)$-strategically closed. Therefore, $\omega_1$ and $\omega_2$ are preserved after forcing with $\mathbb{P}$.

Choose a $V$-generic set $G\subseteq\mathbb{P}$.
For every $\alpha\in S^{\omega_2}_\omega$ let $S_\alpha=p(\alpha)$ for some (any) $p\in G$ with $\alpha\in \mathrm{dom}(p)$.
Observe that $(S_\alpha:\alpha\in S^{\omega_2}_\omega)$ exemplifies $(\ast)$ in $V[G]$.
We claim that $V[G]$ is a model of \textsf{PFA}.

To see this, suppose that $\name{\mathbb{Q}}$ is a $\mathbb{P}$-name of a proper forcing notion and $\name{D}=(\name{D}_i:i\in\omega)$ is a $\mathbb{P}$-name for a family of $\omega_1$ dense subsets of $\name{\mathbb{Q}}$.
Given $p\in\mathbb{P}$ we have to find a condition $t\leq p$ which forces the existence of a filter on $\name{\mathbb{Q}}$ meeting every $\name{D}_i$.
We will define a $\mathbb{P}$-name for an auxiliary forcing notion $\name{\mathbb{R}}$ and we will show that $\mathbb{P}\ast\name{\mathbb{Q}}\ast\name{\mathbb{R}}$ is proper. This $\name{\mathbb{R}}$ will help us define $t$.
 
In $V^{\mathbb{P}}$, a condition $f\in\name{\mathbb{R}}$ is an increasing and continuous function $f:\alpha\rightarrow\omega_2^V$, where $\alpha\in \omega_1$, satisfying that if $\delta\in{\rm dom}(f)$ and $\delta$ is a limit ordinal then $f''\delta\subseteq S_{f(\delta)}$.
The order of $\name{\mathbb{R}}$ is reverse inclusion. 
%
%Since $\mathbb{P}$ is countably closed, $\{(p, \name{r})\in \mathbb{P}\ast \name{R}: \exists r\in V \ p\Vdash_\mathbb{P} \name{r}=\check{r}\}$ is dense in $\mathbb{P}\ast \name{\mathbb{R}}$.
\begin{claim}\label{clmproper}
$\mathbb{P}\ast\name{\mathbb{Q}}\ast\name{\mathbb{R}}$ is proper.
\end{claim}
\par\noindent\emph{Proof}. \newline 
Suppose that $\chi$ is a sufficiently large regular cardinal, $(p,\name{q},\name{r})\in \mathbb{P}\ast\name{\mathbb{Q}}\ast\name{\mathbb{R}}, M\prec\mathcal{H}(\chi)$ is countable and $(p,\name{q},\name{r}), \mathbb{P}\ast\name{\mathbb{Q}}\ast\name{\mathbb{R}}\in M$.
We want to extend $(p,\name{q}, \name{r})$ to an $(M, \mathbb{P}\ast\name{\mathbb{Q}}\ast\name{\mathbb{R}})$-generic condition.

Let $\delta=M\cap\omega_1\in\omega_1$ and let $\gamma=\sup(M\cap\omega_2)$ so $\cf(\gamma)=\omega$.
Since $\mathbb{P}$ is countably closed we can choose an $M$-generic sequence $(p_n:n\in\omega)$ such that $p_0=p$.
Let $p^*=\bigcup_{n\in\omega}p_n\cup\langle\gamma,M\cap\omega_2^V\rangle \in\mathbb{P}$.
Using the fact that $\name{\mathbb{Q}}$ is a $\mathbb{P}$-name of a proper forcing, we can choose a $\mathbb{P}$-name $\name{q}^*$  so that $p^*$ forces that $\name{q}^*$ extends $\name{q}$ and $\name{q}^*$ is $(M[G],\name{\mathbb{Q}})$-generic.

By the choice of $p^*,\name{q}^*$ we see that $(p^*,\name{q}^*)\Vdash M[G\ast H]\cap{\rm Ord}=M\cap{\rm Ord}$. Let $G*H$ be generic for $\mathbb{P}\ast \name{\mathbb{Q}}$ over $V$, containing $(p^*, \name{q}^*)$.
Since $M[G\ast H]$ is countable, we can build a decreasing generic sequence $\bar{r}=(r_i:i\in\omega)$ of conditions in $\mathbb{R}$ over $M[G\ast H]$.
We define $r^*=\bigcup_{n\in\omega}r_i\cup\langle\delta,\gamma\rangle$. We claim that $r^*\in \mathbb{R}$. Indeed, $r^*$ is continuous because $\bar{r}$ is generic over $M[G*H]$ and ${r^*} ''\delta= \bigcup_{i\in\omega} r_i'' \delta \subset M[G*H]\cap \omega_2^V = M\cap \omega_2^V = p^*(\gamma)=S_\gamma=S_{r^*(\delta)}$.

Let $\name{r}^*$ be a $\mathbb{P}*\name{\mathbb{Q}}$-name for a condition in $\name{\mathbb{R}}$ such that $(p^*, \name{q}^*)$ forces that $\name{r}^*$ has the properties as above. Then condition $(p^*,\name{q}^*,\name{r}^*)\leq_{\mathbb{P}\ast\name{\mathbb{Q}}\ast\name{\mathbb{R}}} (p,\name{q}, \name{r})$ is $(M, \mathbb{P}\ast\name{\mathbb{Q}}\ast\name{\mathbb{R}})$-generic.

\hfill \qedref{clmproper}

Fix a condition $p_0\in\mathbb{P}$.
Let $N\prec\mathcal{H}(\chi)$ be of size $\aleph_1$ such that $\omega_1\subseteq N$, so $N\cap\omega_2=\delta \in S^{\omega_2}_{\omega_1}$, and such that $p_0, \name{D}\in N$ and $\mathbb{P}\ast\name{\mathbb{Q}}\ast\name{\mathbb{R}}\in N$.

We apply \textsf{PFA} in the ground model to $\mathbb{P}\ast\name{\mathbb{Q}}\ast\name{\mathbb{R}}$ and obtain a generic filter $F\subseteq \mathbb{P}/p_0\ast\name{\mathbb{Q}}\ast\name{\mathbb{R}}$ over $N$. Define $t=\bigcup\{p\in\mathbb{P}:\exists\name{q}\exists\name{r}, (p,\name{q},\name{r})\in F\}$.
Notice that $t\in\mathbb{P}$ by virtue of the third coordinate of elements in $F$, and $t\leq_{\mathbb{P}}p_0$. More precisely, let $d'_\delta=_{def}\bigcup \{r: \exists (p, \name{q}, \check{r})\in F\}$ and $d_\delta=d'_\delta [\omega_1]$. There is no generality lost here since $\{(p, \name{r})\in \mathbb{P}* \name{\mathbb{R}}: \exists r\in V \ p\Vdash_\mathbb{P} \name{r} =\check{r}\}$ is a dense subset of $\mathbb{P}* \name{\mathbb{R}}$. Then it is easy to see that $d_\delta$ is a club in $\delta$ of order type $\omega_1$ and for each $\alpha\in \lim d_\delta$, it is the case that $d_\delta\cap \alpha \subset t(\alpha)$. Therefore, $\lim d_\delta$ is a witness to the fact that $t\in \mathbb{P}$.

Therefore, $t \Vdash_{\mathbb{P}} \{\name{q} \in \name{\mathbb{Q}} : \exists (p,\name{q},\name{r})\in F\}$ is a generic filter for $\name{\mathbb{Q}}$ meeting all of $(\name{D}_i: i\in \omega_1)$.

\hfill \qedref{thmpfa}

As we will see in Remark \ref{star}, MM implies the failure of $(*)$.

\section{A positive result for $<\omega_1$-type bounded colorings under Martin's Maximum}\label{mmimplies}

Recall that \textsf{MM} implies $\omega_2\rightarrow^{\rm poly}(\omega_1-cl)^2_{<\omega-bdd}$ as shown in \cite{MR2354904}.
We improve this result by relaxing the boundedness assumption.

\begin{theorem}\label{MMtypebounded}
\label{thmtbdmm} $\mathsf{MM}$ implies $\omega_2\rightarrow^{\rm poly}(\omega_1-cl)^2_{<\omega_1-t-bdd}$.
\end{theorem}

\par\noindent\emph{Proof}. \newline 
Suppose that $f:[\omega_2]^2\rightarrow\omega_2$ is $<\omega_1$-type bounded.
Fix an ordinal $\gamma\in\omega_1$ such that ${\rm otp}(t_{\alpha\beta})\leq\gamma$ for every $\alpha<\beta<\omega_2$.
We shall define a stationary-preserving forcing $\mathbb{Q}=\mathbb{Q}_f$ which forces a closed $f$-rainbow set of order type $\omega_1$.

Let $M$ be a countable elementary submodel of $\mathcal{H}(\chi)$, where $\chi$ is a sufficiently large regular cardinal.
We denote $\sup(M\cap\omega_2)$ by $\delta_M$.
Let $K$ be the structure $(\mathcal{H}(\chi),\in,<^*,f,\gamma)$, where $<^*$ is a fixed well-ordering of $\mathcal{H}(\chi)$. Let $\mathcal{S}$ be the collection of countable elementary submodels of $K$. Hence for any $N\in \mathcal{S}$, $\gamma< N\cap \omega_1$.
We define a forcing notion $\mathbb{P}$ as follows.
A condition $p\in\mathbb{P}$ is a pair $(\ell_p,\mathcal{N}_p)$ such that:
\begin{enumerate}
\item [$(a)$] $\mathcal{N}_p=(N_i:i<k), N_i \in \mathcal{S}$ for every $i<k$, $N_i\in N_{i+1}$ for every $i<k-1$.
\item [$(b)$] $\ell_p:\mathcal{N}_p\rightarrow\omega$, and for each $n\in\omega$, $a^p_n=_{def}\{\delta_N:N\in{\rm dom}(\ell_p), \ell_p(N)=n\}$ is  $f$-rainbow.
\end{enumerate}
If $p,q\in\mathbb{P}$ then $q\leq_{\mathbb{P}}p$ iff $\mathcal{N}_p\subseteq\mathcal{N}_q$ and $\ell_q\upharpoonright\mathcal{N}_p=\ell_p$.

We shall use the following terminology for describing a sort of a one-point extension.
If $p\in\mathbb{P}$ then the pair $(M,n) \in \mathcal{S}\times \omega$ is \emph{addable} to $p$ iff $p\in M$ and $a^p_n\cup\{\delta_M\}$ is $f$-rainbow.
If $(M,n)$ is addable to $p$ then $q=p+(M,n)$ will be the condition $(\ell_q,\mathcal{N}_q)$ where $\mathcal{N}_q=\mathcal{N}_p\cup\{M\}, \ell_q\upharpoonright\mathcal{N}_p=\ell_p$ and $\ell_q(M)=n$.

\begin{claim}
\label{mmproperness}
$\mathbb{P}$ is proper and if $G\subset \mathbb{P}$ is generic over $V$, then $S_G=_{def}\{\delta_N: \exists p\in G, N\in \mathcal{N}_p\}$ is a stationary subset of $\omega_2^V$.
\end{claim}
\par\noindent\emph{Proof}. \newline 
The proof is essentially the same as that in \cite{MR716846} and \cite{MR2354904}. But we include a short proof adapted to the present context for completeness. To see that $\mathbb{P}$ is proper, let $\theta>\chi$ be a large enough regular cardinal and let $M^*\prec H(\theta)$ be a countable elementary submodel containing $\mathbb{P}, K$. Given $p\in M^*\cap \mathbb{P}$, since $p$ is finite, there exists $m_0\in \omega$ larger than $\ell_p '' \mathcal{N}_p$. Let $M=M^*\cap K \in \mathcal{S}$. We claim that $q=p+(M, m_0)$ is $(M^*, \mathbb{P})$-generic. To see this, let $r\leq q$ and a dense open set $D\in M^*\cap \mathcal{P}(\mathbb{P})$ be given. We may assume $r\in D$. The goal is to find $t\in D\cap M^*$ such that $t$ is compatible with $r$. Note that $r$ can be decomposed into the part inside $M$, denoted as $r\restriction M$, and the part outside $M$. To simplify presentations, we will assume the part of $r$ outside $M$ consists of $\in$-increasing $\{M_0, M_1, M_2\}$ with $M_0=M$ and $\ell_r(M_i)=m_i$ for $i<3$. Consider $$B=_{def} \{(\alpha_0,\alpha_1,\alpha_2)\in [\omega_1]^3: \exists t\in D, \mathcal{N}_t =(\mathcal{N}_r\cap M) \cup \{M_0', M_1', M_2'\}, $$ $$\ell_t\restriction M =\ell_r\restriction M, \forall i<3 \ \ell_t(M'_i)=m_i, \alpha_i = \delta_{M'_i}\}.$$
Clearly $B$ is definable from $D$ and $r\restriction M$. So $B\in M^*$. As $B\in H(\chi)$ as well, we know that $B\in M$, hence $B\in M_i$ for $i<3$. Let $Qx$ be the Keisler's quantifier asserting ``there exist uncountably many $x$''. Note that $Q\alpha_0 \ Q\alpha_1 \ Q \alpha_2 \ (\alpha_0, \alpha_1, \alpha_2)\in B$ holds. To see this, let $\delta_i=\delta_{M_i}$ for $i<3$, then we know $(\delta_0, \delta_1, \delta_2)\in B$. Work in $M_2$, we know by elementarity that $Q \alpha_2$ $(\delta_0, \delta_1, \alpha_2)\in B$. Repeat the same argument with $M_1$ and $M_0$, we will get $Q\alpha_0 \ Q\alpha_1 \ Q \alpha_2 \ (\alpha_0, \alpha_1, \alpha_2)\in B$.

By elementarity of $M$, $Q\alpha_0 \ Q\alpha_1 \ Q \alpha_2 \ (\alpha_0, \alpha_1, \alpha_2)\in B$ holds in $M$. We say $\nu$ is \emph{good} for $A$ if for any $\{\xi,\zeta\}_<\in [A]^2$, if $\nu<\zeta$ then $\nu \not \in t_{\xi\zeta}=_{def}\{\eta: f(\eta, \zeta)=f(\xi,\zeta)\}$. Consider $C=\{\delta_N: N\in \mathcal{N}_r\}$. Since $C$ is finite and all the color classes have order type bounded above by $\gamma$, the order type of $\bigcup_{\{\xi, \zeta\}\in [C]^2} t_{\xi\zeta}$ is bounded above by $(\gamma+\omega)^\omega$ (this is the ordinal exponentiation), which is $<M\cap \omega_1$. Therefore, we can find $\alpha_0\in M$ that is good for $C$ such that $Q\alpha_1 \ Q \alpha_2 \ (\alpha_0, \alpha_1, \alpha_2)\in B$ holds. Repeat the same argument, we can continue to find $\alpha_1 \in M$ good for $C\cup \{\alpha_0\}$, and $\alpha_2 \in M$ good for $C\cup \{\alpha_0, \alpha_1\}$. Finally, let $t\in D\cap M$ witness $\{\alpha_0, \alpha_1, \alpha_2\}\in B$. It is easy to see that $t$ and $r$ are compatible.

To see $S_G$ is a stationary subset of $\omega_2^V$, let $\name{C}$ be a $\mathbb{P}$-name for a club in $\omega_2^V$ and $p\in \mathbb{P}$ be given. Let $\theta$ be a large enough regular cardinal, and $M^*\prec H(\theta)$ be a countable elementary submodel containing $\mathbb{P}, K, \name{C}, p$. Let $q\leq_\mathbb{P} p$ be an $(M^*, \mathbb{P})$-generic condition, with $M^*\cap K \in \mathcal{N}_q$. Then $q\Vdash_\mathbb{P} \sup (M\cap \omega_2^V)=\sup (M^*\cap \omega_2^V) =\sup (M^*[\name{G}]\cap \omega_2^V)\in \name{C}$.

\hfill \qedref{S_p} \newline

Let $G\subset \mathbb{P}$ be generic over $V$. In $V[G]$, $\omega_2^V$ is collapsed and has cofinality $\omega_1$.
Let $\ell=\bigcup\{\ell_p:p\in G\}$, so $\ell$ partitions $S_G$ into countably many $f$-rainbow subsets. Therefore, at least one of them is stationary in $\omega_2^V$.
Let $\name{n}$ be a name of the least natural number for which the set $S=\ell^{-1}(\{n\})$ is stationary in $\omega_2^V$.
In $V[G]$, let $\mathbb{R}$ be a forcing notion shooting a club of order type $\omega_1$ into $S$. More precisely, $\mathbb{R}$ consists of countable closed subsets of $S$, ordered by end extension.
In $V$, let $\name{\mathbb{R}}$ be the corresponding $\mathbb{P}$-name.  Let $\mathbb{Q}=\mathbb{P}\ast\name{\mathbb{R}}$.

We shall prove that $\mathbb{Q}$ preserves stationary subsets of $\omega_1$.
This will conclude the proof since the generic set for $\mathbb{Q}$ gives a closed $f$-rainbow subset of order type $\omega_1$.
Suppose that $(p,\name{r})\in\mathbb{Q}$ and by extending $p$ if needed we may assume that $p\Vdash_\mathbb{P}\name{n}=n$.
Define $S_p=\{\delta\in\omega_2:\exists M \in \mathcal{S}, \ \delta_M=\delta$ and $(M,n)$ is addable to $p\}$.
Notice that $S_p\subseteq S^{\omega_2}_{\omega}$.

\begin{claim} \label{S_p}
$S_p$ is stationary.
\end{claim}
\par\noindent\emph{Proof}. \newline 
Suppose for the sake of contradiction that $C\subseteq\omega_2$ is a club of $\omega_2$ disjoint from $S_p$. Let $G\subset \mathbb{P}$ containing $p$ be generic over $V$. By an easy density argument, we know for a tail $N\in S_G$, $p\in N$. Let $M\in S_G$ containing $p$ such that $\delta_M\in C\cap \ell^{-1}(n)$. By definition, $(M,n)$ is addable to $p$ and hence $\delta_M\in S_p$. This contradicts with our hypothesis.

\hfill \qedref{mmproperness} \newline

Let $T$ be a stationary subset of $\omega_1$ in the ground model and suppose that $(p,\name{r})\in\mathbb{Q}$.
Let $\name{C}$ be a $\mathbb{Q}$-name for a club in $\omega_1$. Let $\theta$ be a sufficiently large regular cardinal.
Choose $M^* \prec H(\theta)$ containing $(p,\name{r}), \mathbb{Q}, K, \name{C}$ such that letting $M=_{def}M^*\cap K$, we have $(M,n)$ is addable to $p$ and $M\cap\omega_1\in T$. This is possible by Claim \ref{S_p} and the observation that if $M_0, M_1\in \mathcal{S}$ with $p\in M_0\cap M_1$ and $\delta_{M_0} = \delta_{M_1}$, then $(M_0, n)$ is addable to $p$ iff $(M_1, n)$ is addable to $p$.

Let $q=p+(M,n)$.
Notice that $q$ is $(M^*,\mathbb{P})$-generic by Claim \ref{mmproperness} and therefore $q\Vdash_{\mathbb{P}}M^*[G]\cap{\rm Ord}=M^*\cap{\rm Ord}$.
In particular, $\sup(M^*[G]\cap{\omega_2^V})=\sup(M^*\cap{\omega_2^V})\in \ell^{-1}(n)$.
It follows that $q$ forces that there is a decreasing generic sequence over $M^*[G]$, say $\{r_n \in \mathbb{R} :n\in\omega\}$ with $r_0=r$ and a master condition $t$ which satisfies $t\leq_\mathbb{R} r_n$ for every $n\in\omega$.
Since $(q,\name{t})\Vdash_{\mathbb{Q}} M[G\ast H]\cap\omega_1=M\cap\omega_1\in\name{C}\cap T$, we are done.

\hfill \qedref{thmtbdmm}

\begin{remark}\label{star}
Adapting the proof of Theorem \ref{MMtypebounded}, one can also show  the following: for any $f: [\omega_2]^2\to \omega_2$ such that there exists a stationary $S\subset S^2_0$ and $\gamma<\omega_1$ with the order type of $t_{\eta\alpha}=\{\eta'<\alpha: f(\eta', \alpha)=f(\eta, \alpha)\}$ being bounded above by $\gamma$ for any $\alpha\in S$ and $\eta<\alpha$, then there exists a closed $f$-rainbow set $C\subset S$ of order type $\omega_1$. In particular, this shows that MM implies $\neg (*)$. To see this, suppose a putative $(*)$-sequence $\langle S_\alpha: \alpha\in S^2_0\rangle$ is given. We may find a stationary $S\subset S^2_0$ and $\gamma<\omega_1$ such that for all $\alpha\in S$, $\mathrm{otp}(S_\alpha)\leq \gamma$. Define a  coloring $f$ on $[\omega_2]^2$ such that for each $\alpha\in S$, 
\begin{itemize}
\item $f(\cdot, \alpha)\restriction S_\alpha$ is constant and 
\item for any $\alpha\in S$ and any $\eta<\alpha$, $\mathrm{otp}(t_{\eta\alpha})\leq \gamma$.
\end{itemize}
Thus $f$ satisfies the modified hypothesis above as witnessed by $S$ and $\gamma$. Therefore, MM implies there is a closed $f$-rainbow $C\subset S$ of order type $\omega_1$. Let $\beta=\sup C$ and $c_\beta$ be the witness to the property of the $(*)$-sequence (recall Definition \ref{clmstar}). Letting $D=C\cap c_\beta$, we see that on one hand, for any $\alpha\in D$, $D\cap \alpha\subset S_\alpha$, but on the other hand, for any $\alpha\in D$, $|D\cap S_\alpha|\leq 1$. We have reached a contradiction.
\end{remark}

%
%There are constraints for Theorem \ref{MMtypebounded} at larger cardinals. 
%
%\begin{claim}
%\label{clmomega3} If $\kappa\geq\omega_2$ then $\kappa\nrightarrow^{\rm poly}(\omega_1+1-cl)^2_{\kappa,\omega_1-t-bdd}$ and even $\kappa\nrightarrow^{\rm poly}(\omega_1+1-st)^2_{\kappa,\omega_1-t-bdd}$.
%\end{claim}
%
%\par\noindent\emph{Proof}. \newline 
%Let $S=\{\delta\in\kappa:\cf(\delta)=\omega_1\}$.
%We define this set even if $\cf(\kappa)=\omega$.
%For every $\beta\in S$ fix a club $c_\beta$ of $\beta$ so that ${\rm otp}(\beta)=\omega_1$.
%Let $A=\{j_\beta:\beta\in S\}\subseteq\kappa$ be such that $|\kappa-A|=\kappa$.
%One can think of $A$ as a set of colors, and we leave enough room for $\kappa$ many additional colors.
%
%For every $\beta\in S$ and every $\alpha\in c_\beta$ we let $f(\alpha,\beta)=j_\beta$.
%For every other pair $(\gamma,\delta)\in[\kappa]^2$ let $f(\gamma,\delta)$ be a new color from $\kappa-A$.
%Clearly $f$ is $\omega_1$-type bounded.
%We claim that $f$ exemplifies the negative asserted relation.
%
%To see this, suppose that $C\subseteq\kappa,{\rm otp}(C)=\omega+1$.
%Let $\beta$ be the top element of $C$.
%Since both $C\cap\beta$ and $c_\beta$ are club subsets of $\beta$ and $\cf(\beta)>\omega$ we can choose $\alpha,\alpha'\in C\cap c_\beta$ such that $\alpha\neq\alpha'$.
%By definition, $f(\alpha,\beta)=f(\alpha',\beta)=j_\beta$, so $C$ fails to be $f$-rainbow.
%The same argument works if $C\cap\beta$ is a stationary subset of $\beta$, so we are done.
%
%\hfill \qedref{clmomega3}

\section{Some open problems}\label{openquestion}

We have seen that $\omega_2\nrightarrow^{\rm poly}(\omega_2-cl)^2_{2-bdd}$ in \textsf{ZFC}.
On the other hand, it is consistent to get a closed rainbow subset of type strictly greater than $\omega_1$.

\begin{question}
\label{qalphamm} 
Is it consistent that $\omega_2\rightarrow^{\rm poly}(\alpha-cl)^2_{2-bdd}$ for every ordinal $\alpha\in\omega_2$?
\end{question}

Note $\alpha=\omega_1+\omega_1$ is the first point which challenges our method.

\begin{question}
\label{qomega2stat} Is it consistent that $\omega_2\rightarrow^{\rm poly}(\omega_2-st)^2_{2-bdd}$?
\end{question}

%As we have seen, $\omega_2\nrightarrow^{\rm poly}(\omega_2-cl)^2_{\omega_2,2-bdd}$.
%But maybe the correct generalization of the closed rainbow set of size $\omega_1$ obtained under \textsf{MM}, is when the domain of the coloring is bigger than the size of the rainbow set.
%This can be stretched a bit more, as reflected in the following:
%
%\begin{question}
%\label{qalmost} Is it consistent that for some (large) cardinal $\kappa$ we will have $\kappa\rightarrow^{\rm poly}(\alpha-cl)^2_{\kappa,2-bdd}$ for every $\alpha\in\kappa$? \newline 
%Is it consistent that $\kappa\rightarrow^{\rm poly}(\omega_2-cl)^2_{\kappa,2-bdd}$?
%\end{question}
%
%We consider now a wider class of colorings by relaxing the boundedness property.
%We saw that Chang's conjecture is helpful for obtaining positive rainbow relations with infinite degree of boundedness.
One may wonder if Corollary \ref{thmsing} can be improved:

\begin{question}
\label{qthmsing} Is it consistent that $\aleph_{\omega+1}\rightarrow^{\rm poly}(\omega_1-cl)^2_{\aleph_{\omega+1},(<\aleph_\omega)-bdd}$?
\end{question}

We have only dealt with successors of singular cardinals of countable cofinality. Our method is limited because the consistency of the higher analogue of Theorem \ref{pfatypebounded} is not clear (see Question \ref{qomega2stat}).

\begin{question}
Is it consistent that $\kappa$ is a singular cardinal with $cf(\kappa)>\omega$ and $\kappa^+\to (cf(\kappa)^+)^2_{<\kappa-bdd}$?
\end{question}

We saw that $\omega_2\rightarrow^{\rm poly}(\omega_1-cl)^2_{<\omega_1-t-bdd}$ and $\omega_2\rightarrow^{\rm poly}(\omega_1-st)^2_{\omega-bdd}$ are respectively consistent.
A natural question is whether the following joint strengthening is possible.

\begin{question}
\label{qmmctblbd} Is it consistent that $\omega_2\rightarrow^{\rm poly}(\omega_1-cl)^2_{\omega-bdd}$ under any assumption?
\end{question}

\bibliographystyle{plain}
\bibliography{bib}

\end{document}